# Study on the impact of uncertain design parameters on the perfomances of a permanent magnet assisted synchronous reluctance motor

Adán Reyes Reyes[1,2,*], André Nasr[1], Delphine Sinoquet[1], and Sami Hlioui[2,3]

[1] IFP Energies nouvelles, Institut Carnot IFPEN Transports Energie, 1 et 4 avenue de Bois-Préau, 92852 Rueil-Malmaison, France
[2] Paris Saclay University, CNRS, SATIE, Gif-sur-Yvette, France
[3] CY Cergy Paris University, CNRS, SATIE, Cergy, France
* adan.reyes-reyes@ifpen.fr

*Abstract—* In this paper, deterministic and robust design optimizations of a permanent magnet assisted synchronous reluctance machine were performed to study the impact of different uncertain input parameters on the design. These optimizations were carried out using a surrogate model based on 2-D finite element simulations. Different robust optimizations considering geometric and magnetic uncertain parameters were compared to the deterministic optimization. It was noticed that both geometrical and magnetic properties tolerances greatly impact the machines' performances, where the magnetic properties had a more significant impact on the mean torque. In such a case, robust optimization is essential to find optimal and robust electric motors designs.

## 1 Introduction

With the increasing concerns over climate change, many measures have been adopted to reduce greenhouse gas emissions. For transportation systems, in order to replace internal combustion engine vehicles, electric and hybrid vehicles (EV, HEV) have been intensively developed. In these vehicles, the electrical machine is one of their main components.

Among the different types of electrical machines used in electric vehicles, Permanent Magnet assisted Synchronous Reluctance Machines (PMaSynRel) are one of the most used machines nowadays thanks to their good performances and their relatively low cost [1]. Unlike Surface-Mounted Permanent Magnet Synchronous Machines (SMPMaSynRel), PMaSynRel exploit two types of torque: the hybrid torque generated using permanent magnets and, the reluctance torque which makes profit of the machine's saliency. Since SMPMSMs only generate hybrid torque, they need more permanent magnets to achieve the same torque density and tend therefore to be more expensive. For all these machines, optimization is often used in order to find the best design respecting all the required specifications.

Many optimization methodologies applied on electrical machines can be found in literature [2][3]. Most of these methodologies can be described as deterministic since they do not consider any uncertainties on the input parameters. However, in practice, there are many discrepancies between the theoretical and real (measured) values of these parameters. These differences can be caused by manufacturing and assembly tolerances in the prototype as well as by the lack of precision in magnetic properties of used materials. These variabilities can lead to degraded performances compared to the nominal performances simulated in the design phase. To reduce such deviations, the parameter uncertainties should be considered in the optimization procedure.

In opposition to a deterministic optimization, a robust optimization considers two types of input parameters: certain parameters also known as controllable parameters and uncertain parameters. Controllable parameters are the same ones used in a deterministic optimization whereas uncertain parameters are specific to robust optimization techniques. This type of parameter can take varying values due to the associated uncertainties: it is then modeled by a random variable and an associated Probability Distribution Function (PDF). The presence of random input variables for the simulator leads to random output variables and then, random objective and constraint functions. Various formulations of the resulting optimization problem are proposed in the literature [4][5][6] based on expectation, probability, or quantiles of these random variables.

Reliability Based Design Optimization (RBDO) is a method used to obtain optimal and safe designs in the sense that the outputs of certain functions are inside a security domain, described by constraints. A robust or reliable design has therefore a high probability to respect these constraints. Examples of this approach can be found in [7].

Worst-case optimization considers the extreme values as objective functions and/or constraints i.e., the maximum or minimum value of the outputs caused by the uncertainty propagation [8].

There is another very common formulation which was also adopted in this work: Robust Design Optimization (RDO). In this methodology, the expectation (average) of the objective function is optimized. To limit extreme values, a second objective based on the variances of the objective function can be added [9].

The computation of robustness metrics such as expectations or variances requires a large sample of the uncertain input variables and thus a large number of simulations. To limit this high computational cost especially with the use of finite element simulations, meta-modeling techniques coupled with design of experiments are used to replace the costly simulations by predictions using the resulting surrogate models [10].

We will present in this paper a study on the impact of uncertain design parameters on the performances of a PMaSynRel motor. It is a 3-phases 10-poles 60-slots PMaSynRel with a Machaon rotor structure (Fig. 1). It has an outer stator diameter of 220 mm and an active length of 200 mm. Each pole has 3 flux barriers and 7 PMs. This machine was initially designed for an EV application having a maximum torque of 430 N.m. The design of this motor was originally obtained via a deterministic optimization. In order to study the sensitivity of this design and the possibility of improving its robustness as well as its performances, new deterministic and robust optimizations will be performed. To do so, techniques like Design Of Experiments (DOE), Finite Element Method (FEM) surrogate modeling, sensitivity analysis, quasi Monte Carlo methods and optimization algorithms were used.

The remainder of this paper is as follows: We will first introduce the surrogate models-based methodology used in this work. Next, we will analyze the obtained results via deterministic and robust formulation optimizations. Finally, we will verify some solutions' performances obtained by the metamodels through Finite Element simulations.

TABLE I: Optimization variables

| Input Parameter | Lower bound $x_l$ | Upper bound $x_u$ | Manufacturing Tolerance |
|---|---|---|---|
| Slot_angle | 2.47° | 3.27° | ±0.1° |
| Beta_L1_P1 | 27.03° | 29.66° | ±0.33° |
| Beta_L1_P2 | 37.03° | 39.66° | ±0.33° |
| Beta_L2_P1 | 31.03° | 33.66° | ±0.33° |
| Beta_L2_P2 | 47.03° | 49.66° | ±0.33° |
| Beta_L3_P1 | 33.7° | 37° | ±0.33° |
| Beta_L3_P2 | 59.7° | 63° | ±0.33° |
| Airgap | 0.55 mm | 0.65 mm | ±0.03 mm |
| Bridge_L1 | 2.6 mm | 2.98 mm | ±0.05 mm |
| Bridge_L2 | 0.9 mm | 1.18 mm | ±0.05 mm |
| Bridge_L3 | 0.5 mm | 0.62 mm | ±0.03 mm |
| Bridge_tang | 0.4 mm | 0.6 mm | ±0.05 mm |

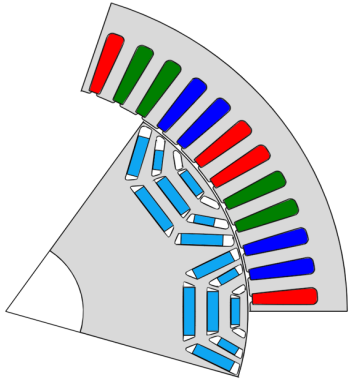

Fig. 1: Geometry of the PMaSynRel motor

## 2 OPTIMIZATION WORKFLOW

As mentioned before, we will perform in this paper deterministic and robust optimizations in order to study the sensitivity of a PMaSynRel machine to design parameters uncertainties. Fig. 2 shows the design parameters for the stator and for the rotor considering one layer. TABLE I lists all the design parameters to be used in the optimizations as well as their lower and upper bounds. The manufacturing tolerance is also given.

Geometrical parameters are not the only design parameters that can be uncertain. The characteristics of the magnetic materials used to build a machine, like permanent magnets and electrical steel, can also present some uncertainties and their characteristics can deviate from those given in the datasheet.

The permanent magnet model used in our electromagnetic simulations is a linear model defined by two quantities: $B_r$ which is the remanent induction and $H_{cb}$ which is the coercitive magnetic field (Fig. 3). We will consider that the values of $B_r$ and $H_{cb}$ can vary between their nominal characteristics and, in the worst case, lower characteristics representing a 6.5 % degradation. A coefficient β between 0% and 6.5% will be used as a random variable describing the level of degradation in the robust optimizations.

As for the B(H) characteristic, the degradation will be limited to its "knee" part. A random variable α between 0 (representing fully degraded characteristic) and 1 (nominal characteristic) will be used to define the degradation level of the B(H) curve. At maximum degradation (α = 0), the induction level is reduced by 25 %. Fig. 4 shows a comparison between the initial and degraded B(H) curves.

Fig. 5 shows the workflow used to carry out the optimizations. At first, a DOE was built with the upper and lower bounds of the input parameters shown in TABLE I to fit a surrogate model for each of the considered objective functions: the mean torque and torque ripple. Secondly, these models were used to perform a global sensitivity analysis to detect the most impacting parameters on the objective functions. This will allow us to limit the number of parameters considered as uncertain. At last, and after performing the meta-model-based deterministic and robust optimizations, FEM simulations were carried out to verify the results.

These 5 parameters will be then considered as uncertain variables in the robust optimization. The dispersion of these variables will be then integrated in the robust optimizations by considering a perturbation vector *U*.

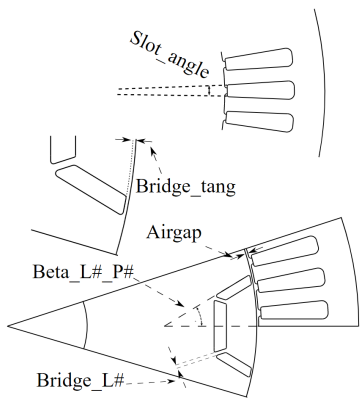

Fig. 2: Design parameters for one layer (# is the number of the layer).

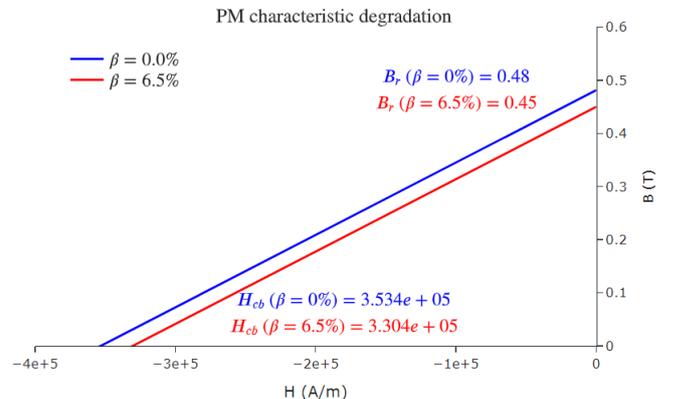

Fig. 3: Nominal and degraded PM characteristic.

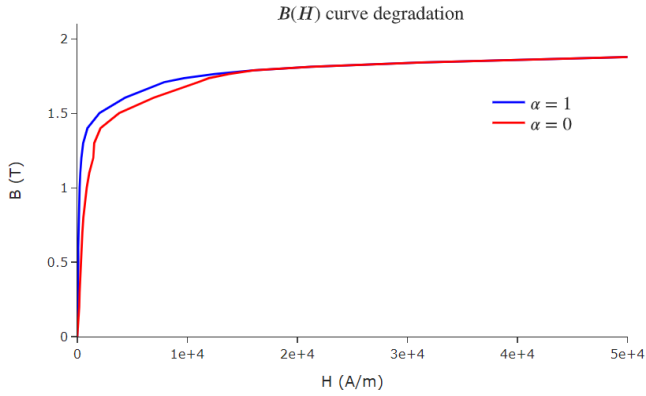

Fig. 4: Nominal and degraded B(H) curve.

## 2.1 Surrogate Models

To reduce computation time, surrogate models have been built for each of the objective functions. To build such models, there are three steps to follow: build a DOE, train the metamodel and check its predictivity with a test set. The chosen DOE is a maximin Latin Hypercube Sampling (LHS) which is a technique that covers well the search space while preserving good projection properties [11]. This DOE was built with 234 points using the bounds described in TABLE I. As for the surrogate model, it is a universal Kriging with linear trend function. For mean torque, a tensorized Matérn 5/2 covariance function has given the best predictivity. As for Torque ripple, a tensorized absolute value exponential kernel was used since the Torque ripple function is not that smooth. Kriging was chosen as metamodel since it is very good at learning nonlinear objective functions and has demonstrated good performances in electrical machines optimization [12]. Finally, to evaluate the accuracy of the metamodel, a Normalized Root Mean Square Error (NRMSE) was computed on a validation test set:

$$NRSME = \frac{||y_{real} - y_{pred}||}{||y_{real}||} * 100\% \quad (1)$$

zero indicates a good model fit. The obtained *NRSME* of the kriging model for mean torque is *0.2%* and for torque ripple is 8%. These results were obtained with a train and test sets composed by 175 and 59 samples, respectively. Torque ripple depends not only on mean torque but also on torque amplitude which makes this function more difficult to model than mean torque. We consider those metamodels sufficiently accurate for performing the sensitivity analysis and the optimization procedures.

## 2.2 Sensitivity Analysis

For the robust optimizations, we have decided not to consider all the geometric parameters given before as uncertain parameters. Only the most impacting ones will be considered. To do so, a sensitivity analysis will be performed [13]. For this work, the Sobol Indices were chosen as they measure the global impact of the input variables on the output functions. Such indices represent the amount (or percentage)

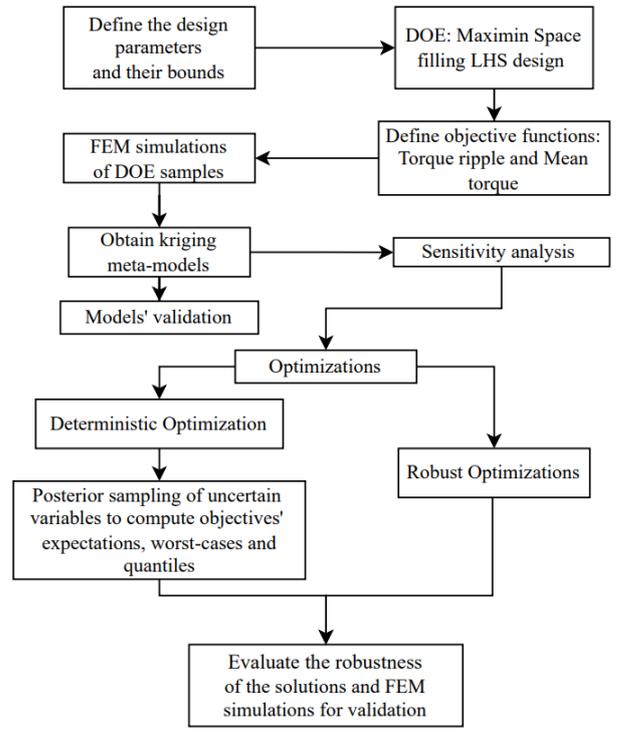

Fig. 5: Optimization workflow.

of the total output variance attributable to each subset of input variables. For instance, in the simplest case of one output variable *Y* and two inputs $X_1$ and $X_2$, we find three factors causing the total variance of *Y*: due to the variation of $X_1$ alone, due to the variation of $X_2$ alone, and due to the variation of $X_1$ and $X_2$ simultaneously. Dividing these three quantities by the total variance of *Y*, we can obtain three percentages that can be directly considered as sensitivity measures. The largest is the value of these indices, the largest is their importance. The commonly used ones are the first order indices ($S_i$) representing the contribution to the variance of *Y* due to the variation uniquely of $X_i$, and the total indices ($S_{TOTAL,i}$) representing the contribution to the variance of Y due to the variation of $X_i$ and all its interactions with the remaining input variables:

$$S_i = \frac{VAR_{X_i}(E_{X_{\sim i}}[Y|X_i])}{VAR(Y)} \quad (2)$$

$$S_{TOTAL\_i} = \frac{E_{X_{\sim i}}[VAR_{X_i}(Y|X_{\sim i})]}{VAR(Y)} \quad (3)$$

where $X_{\sim i} = X_1, \cdots, X_{i-1}, X_{i+1}, \cdots, X_{N_x}$ and $N_x$ is the number of optimization parameters. We computed the indices using the Kriging surrogate models. The results of the sensitivity analysis applied to mean torque and to torque ripple are presented in Fig. 6 and Fig. 7 . Only the most important geometrical inputs are displayed for better visibility. It was found that the stator slot width opening angle (Slot_angle) and the flux barrier opening angles for barriers 1 and 2 (Beta_L1_P1, Beta_L1_P2, Beta_L2_P1 and Beta_L2_P2) have the biggest impacts on the mean torque as well as on the torque ripple. We can see in Fig. 6 that the total order Sobol indices for torque ripple are higher than the first order ones. This means that torque ripple is more sensitive to the interaction between the variables (Beta_L1_P1 and Beta_L1_P2 for example) than to one variable alone. This is

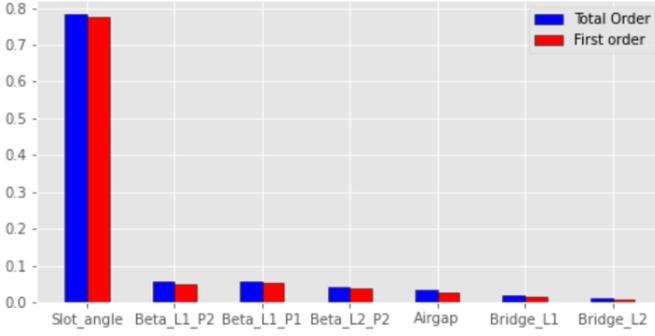
Fig. 6: Mean torque Sobol' indices

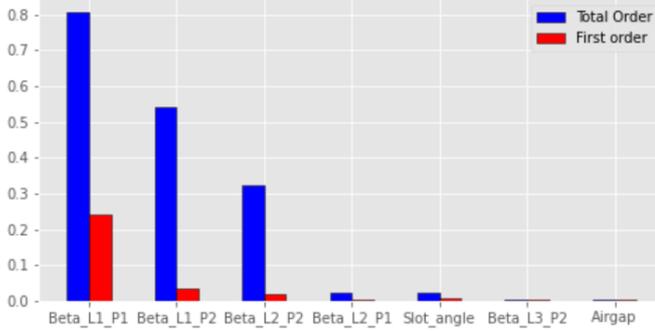
Fig. 7: Torque ripple Sobol' indices

not the case for the mean torque, where the predominant indices are first-order ones.

*2.3 Robustness Metrics*

To perform a robust optimization, one or several robustness measures can be used as objectives or constraints. In this section, we will recall some of them. Consider Fig. 8. This figure shows a possible probability density function (PDF) of the torque ripple considering perturbations on one (or some) input design parameters. The alpha-quantile $q_\alpha$ is the output value for which a given sample has a probability equal to $\alpha$ of being less than or equal to $q_\alpha$. For example, in Fig. 8, the 10% quantile ($q_{10\%}$) has a torque ripple value around 5.5 %. This means that 10 % of the considered samples have a torque ripple less than or equal to 5.5 %. The expectation (E) is equal to the average value of the output variable and would coincide with the 50% quantile in a perfectly symmetric PDF. We can quantify how spread out the output values are by means of the standard deviation (STD), which measures the mean square dispersion around the expectation.

## 3 DETERMINISTIC AND ROBUST OPTIMIZATIONS

We will present in this section the results of different optimizations: a deterministic and a robust optimization. For the deterministic optimization problem, we have:

$$\min_{x \in X} [f_1(x), f_2(x)] \quad (4)$$

where $f_1$ is the opposite of the mean torque (in order to maximize it) and $f_2$ the torque ripple; $X$ is the controllable parameters space defined in TABLE I. For the robust optimization problems, two formulations have been considered:

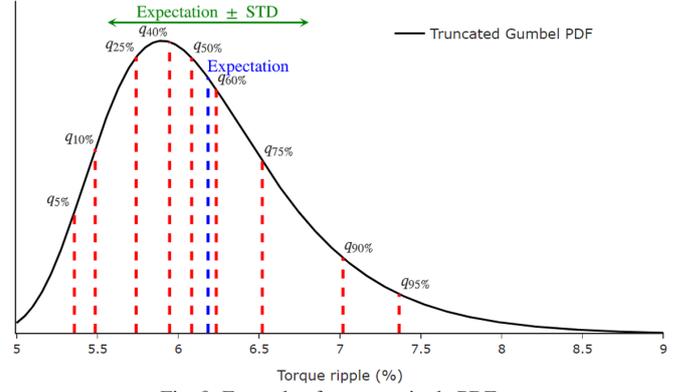
Fig. 8: Example of a torque ripple PDF

- Expectations optimization:

$$\min_{x \in X(U)} [\ E_U[f_1(x+U)], E_U[f_2(x+U)]\ ] \quad (5)$$

- Worst-case optimization:

$$\min_{x \in X(U)} [\ \max_{u \in \Omega} f_1(x+u), \max_{u \in \Omega} f_2(x+u) \quad (6)$$

where $X(U) = [x_{l1}-u_{u1}, x_{u1}+u_{u1}] \times \cdots \times [x_{lNx}-u_{uNx}, x_{uNx}+u_{uNx}]$

and $\Omega = [-u_{u1}, u_{u1}] \times \cdots \times [-u_{uNx}, u_{uNx}] \quad (7)$

$u_{uj}$ is the manufacturing tolerance of the parameter number $j$, and $x_{lj}$ and $x_{uj}$ are the lower and upper bounds for the optimization variable $x_j$, respectively. The tolerance for each geometrical parameter is given in TABLE I. Based on the sensitivity study in the previous section, only 5 geometrical parameters will be considered as uncertain. Their uncertainties were considered following uniform distributions, i.e., $U_j \sim Unif(-u_{uj}, u_{uj})$. For parameters with no considered uncertainties, $u_u$ is simply equal to 0.

The goal of the first robust formulation is to optimize the mean torque's and torque ripple's expectations in a Pareto sense. In the worst-case formulation, the objective is to limit the worst possible value of the mean torque and torque ripple caused by uncertainties (equivalent to $q_{100\%}$ and $q_{0\%}$ respectively). By adopting such formulation, the designer can be sure that all the manufactured machines will exceed the performances found on the Pareto front. On the other hand, a formulation using expectations doesn't guarantee that.

To solve these optimization problems, the genetic algorithm NSGA 2 was used [14]. This algorithm has shown good performances for other studies of electrical machine optimization as in [15]. We used a DOE maximin LHS to compute samples of $x+U$ to calculate the objective functions' expectations with a quasi-Monte Carlo method. When it comes to the Worst-case formulation, we have two options: computing samples of $x+U$ and take the maximum value of these samples as an estimator of $\max_{u \in \Omega} f(x+u)$ *or* obtaining the absolute maximum value with an optimization algorithm. In this work, the last solution was adopted using a Particle Swarm Optimization (PSO) [16] algorithm for the embedded mono-objective optimizations. By doing so, we will obtain more accurate worst-case estimates.

To study the impact of the two different types of uncertainties (related to geometrical parameters and material properties) on the performances of the PMaSynRel machine, we will perform in the following sections two groups of optimizations: the first one will only consider geometrical parameters as uncertain parameters ($U_g$) while the second one will also consider the magnetic material characteristics as uncertain as described before ($U_g$, $U_m$).

### 3.1 Robust optimizations considering uncertainties on geometrical parameters

Fig. 9 shows a comparison between the deterministic (blue) and the Expectations optimization (red) Pareto fronts. The expected performances of the deterministic Pareto front have been reevaluated (pink): The design variables were perturbed by adding sampled values of the uncertain variables. These expected values represent the average mean torque and average torque ripple for each machine obtained by the deterministic Pareto optimization considering a posteriori uncertainty on the input parameters. As we can notice, a deterministic optimization does not guarantee a robust design: For the same average mean torque, the average torque ripple of a batch of machines issued from the deterministic optimization (pink) shows higher values than that of a batch of machines issued from the robust optimization (red). We can also note that the obtained results outperform the initial machine's performances as shown in TABLE II which lists the initial machine's nominal performances as well as its expectations and worst-case values when considering uncertain parameters.

The pareto front of the Worst-case optimization (green) is presented in Fig. 10. As in Fig. 9, the deterministic pareto front was also added (blue). The worst-case performances of the deterministic Pareto front have been evaluated (light green) thanks to a posteriori uncertainty on the input parameters and PSO maximization. Once again, these results show the importance of a robust optimization in limiting the performance degradation that a sample of machines can have.

To go deeper into this analysis, we empirically compared the distribution of different designs. For this purpose, we show in Fig. 11 and Fig. 12 boxplots of a subset of points selected from the deterministic and robust Pareto fronts shown in Fig. 9 and Fig. 10, respectively. Each pair of boxplots in Fig. 11 represents a comparison of the distribution of torque ripple values between a deterministic machine (pink in Fig. 9) and a robust machine (red in Fig. 9) falling in one of the zones ($A_{Exp}^{Ug}$, $B_{Exp}^{Ug}$, $C_{Exp}^{Ug}$, $D_{Exp}^{Ug}$, $E_{Exp}^{Ug}$). The "Exp" index means that the design is issued from a robust optimization based on the expectance formulation. The "Ug" is to say that the considered uncertainties are geometrical parameters. The same boxplots are represented in Fig. 12. However, for this time, the designs were selected from the deterministic Pareto front in Fig. 10 (light green) and the robust one (dark green). All the boxplots in Fig. 11 and Fig. 12 show the values of $q_{25\%}$ (q1), $q_{50\%}$, which is the median of the sample (q2) and $q_{75\%}$ (q3) as represented in Fig. 11 (boxplot zone $C_{Exp}^{Ug}$).

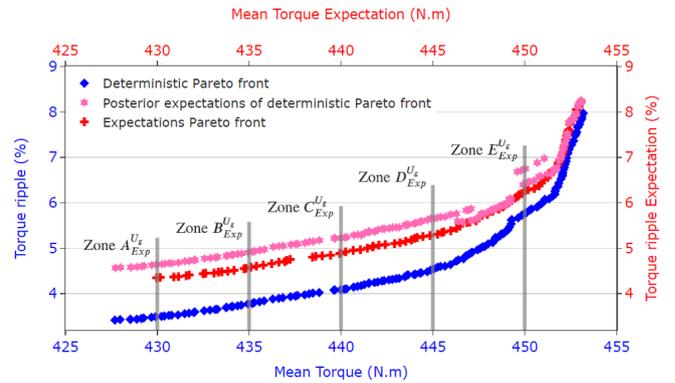

Fig. 9: Pareto fronts obtained by deterministic (front in blue and expectations obtained by posterior perturbations in pink), and expectations optimization (red) where uncertainties come from geometrical perturbations $U_g$. Dark gray zones highlight points with similar Mean Torque expectation values.

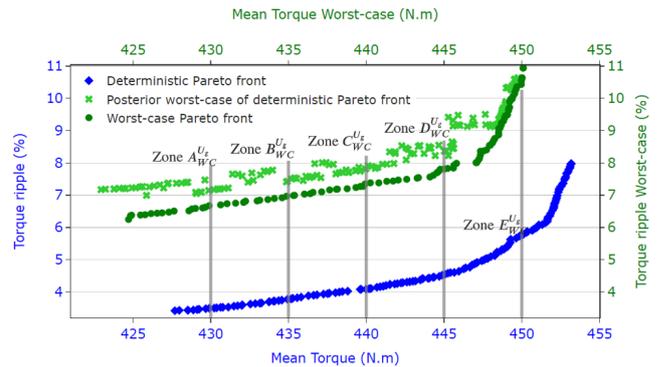

Fig. 10: Pareto fronts obtained by deterministic (front in blue and worst-cases obtained by posterior perturbations in light green), and worst-cases optimization (green) where uncertainties come from geometrical perturbations $U_g$. Dark gray zones highlight points with similar Mean Torque worst-case values.

TABLE II : Initial machine performances

|  | Mean Torque | Torque ripple |
| --- | --- | --- |
| Nominal values | 433.34 (N.m) | 10.38 (%) |
| Expectation ($U_g$) | 432.32 (N.m) | 9.4 (%) |
| Worst-case ($U_g$) | 425.89 (N.m) | 11.41 (%) |
| Expectation ($U_g$,$U_m$) | 422.21 (N.m) | 9.16 (%) |
| Worst-case ($U_g$,$U_m$) | 405.65 (N.m) | 11.42 (%) |

We can observe in Fig. 11 that for each pair of machines, the robust one shows better overall performance than the deterministic one. For example, for the machines in zone $A_{Exp}^{Ug}$, we can notice that the median of the deterministic machine is almost the same as the q3 value of the robust machine (4.5%). This means that there is a 50% chance for a deterministic design to have a torque ripple value higher than 4.5 %. For the robust design, the probability is only 25 %. We can also notice that the expectation and Standard Deviation (STD) of the torque ripple associated with the robust optimization solution (4.3%, 0.35%, respectively) outperform the expectation and STD obtained by the posterior analysis of the deterministic solution (4.6%, 0.56%, respectively). These results stress the importance of robust optimizations when dealing with uncertainties.

Fig. 12 also shows pairs of boxplots comparing predicted torque ripple values between deterministic (light green) and robust (dark green) machines from Fig. 10 with the same minimum (worst-case) predicted mean torque value in each zone. For all the selected zones, the worst-case torque ripple value (maximum value for a batch) of a robust design is lower (or equal) than a deterministic one. Besides this, robust solutions also have much lower STD values than deterministic solutions, especially for zones with low worst-case mean torque. In zone $A_{WC}^{U_g}$ for example, the robust machine (first green boxplot in Fig. 12) has a STD of 0.26% and a worst-case torque ripple of 6.7% compared to 0.54% and 7.2%, respectively for its deterministic counterpart (first blue boxplot in Fig. 12). Both optimizations lead practically to the same design for high values of worst-case mean torque like in zone $E_{WC}^{U_g}$. This can also be seen in Fig. 11 and with the Pareto fronts of the deterministic and robust optimizations getting very close with increasing torque.

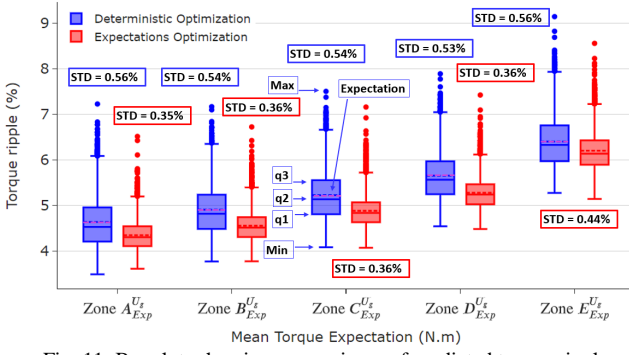

Fig. 11: Boxplots showing comparisons of predicted torque ripple values between deterministic (blue) and robust (red) machines from **Fig. 9** with similar predicted Mean Torque Expectation values.

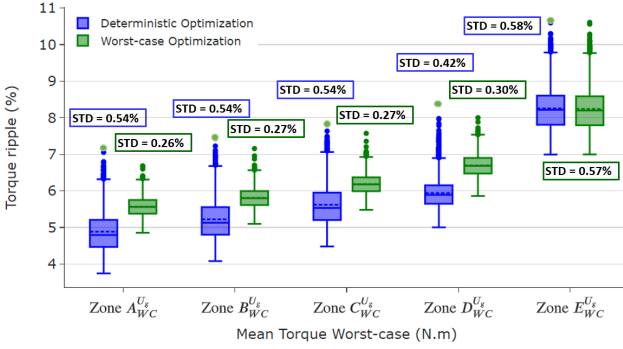

Fig. 12: Boxplots showing comparisons of predicted torque ripple values between deterministic (blue) and robust (green) machines from Fig. 10 with similar predicted Mean Torque Worst-case values.

The worst-case optimization has allowed to limit the performance degradation of the least performant machine in a sample, it has led to worse expectance values of torque ripple than those obtained by deterministic optimization. Instead of using the worst-case as an objective function, it could be used as a constraint in a constrained optimization problem while still using the expectations as an objective. Such formulation allows having good machines samples (expectation-wise) while limiting the worst performances we can have. Given that the quantile and worst-case are related (as the worst-case is equal to the 100% quantile in a minimization scheme), and we have found that the worst-case value could be set as a constraint; a quantile formulation (like 95% quantile) could be used as a more permissive constraint. This is important because finding samples having 100% chance of satisfying a constraint could be difficult, especially when the tail of the PDF of the quantity of interest is large.

### 3.2  Robust optimizations considering uncertainties on geometrical parameters and material properties

In order to evaluate the impact of magnetic material uncertainties on the performances of the PMaSynRel motor, Fig. 13 and Fig. 14 show the results of two robust optimizations performed using respectively the expectations and worst-case formulations. This time, the uncertain parameters are both the geometric parameters and the magnetic material properties. We added on the same figures the Pareto fronts from Fig. 9 and Fig. 10 for comparison.

Compared to results obtained considering only $U_g$ as random variables, we can notice that the introduction of magnetic material uncertainties in the robust optimizations only affect the mean torque while the torque ripple stays almost the same. We can observe that the main effect of degradation of the materials' magnetic properties, is a horizontal translation towards lower values of the mean torque compared to the solutions obtained when we only considered $U_g$ as random variables. In order to look deeper into this result, we present in Fig. 15 the evolution of the estimated mean torque (metamodel) for machine $A_{Exp}^{Ug}$ in Fig. 9 in respect to α and β. As it can be seen, the mean torque metamodel has a linear behavior with respect to these two properties and its maximum value is found when no degradation is applied (α = 1 and β = 0). At α = 0 and β = 1 (full degradation), the torque is around 411 N.m with a 20 N.m decrease.

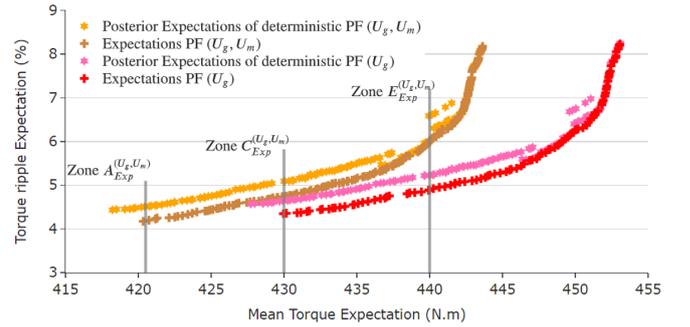

Fig. 13: Pareto fronts obtained by expectations optimization (brown) and expectations obtained by posterior perturbations from deterministic optimization results (orange) where uncertainties come from geometrical and materials' properties perturbations ($U_g$, $U_m$). Red and pink sets of points come from Fig. 9. Dark gray zones highlight points with similar Mean Torque expectation values.

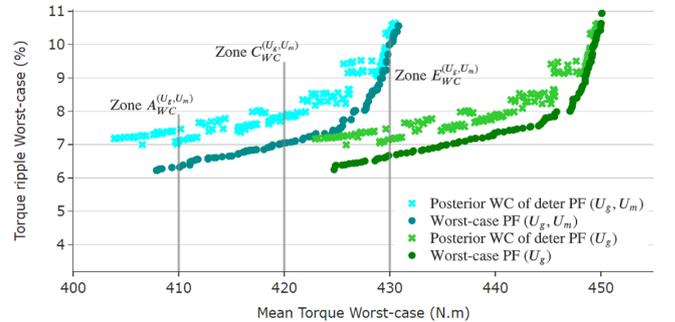

Fig. 14: Pareto fronts obtained by worst-cases optimization (darkcyan) and worst-case obtained by posterior perturbations from deterministic optimization (cyan) where uncertainties come from geometrical and materials' properties perturbations ($U_g$, $U_m$). Green and light green sets of points come from Fig. 10. Dark gray zones highlight points with similar Mean Torque worst-case values.

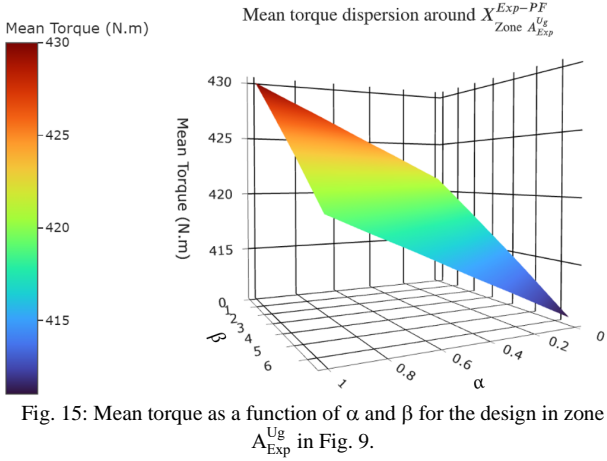

Fig. 15: Mean torque as a function of α and β for the design in zone $A_{Exp}^{Ug}$ in Fig. 9.

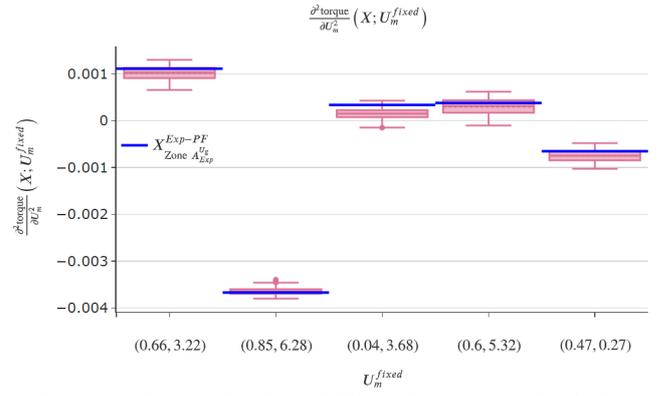

Fig. 16: Boxplots showing the variability (with respect to a batch of Xg samples) of the mean torque metamodel second partial derivative (with respect to Um) evaluated at a given value of Um. The blue horizontal lines stand for the machine related to zone $A_{Exp}^{Ug}$ and coming from the expectations Pareto front.

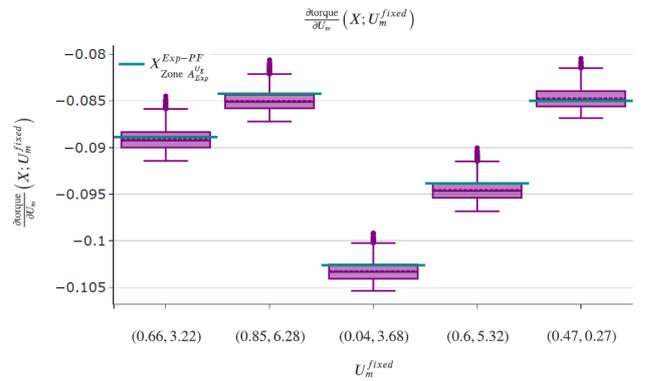

Fig. 17: Boxplots showing the variability (with respect to a batch of Xg samples) of the mean torque metamodel first partial derivative (with respect to Um) evaluated at a given value of Um. The dark cyan horizontal lines stand for the machine related to zone $A_{Exp}^{Ug}$ and coming from the expectations Pareto front.

To confirm this observation all over the search space ($X_g$), a subset of 5000 points was analyzed (this subset was created using LHS maximin over $X_g$ including as a reference to the machine related to zone $A_{Exp}^{Ug}$ and coming from the expectations Pareto front, for which the relationship is linear (as seen in Fig. 15). For all these points, the second partial derivative with respect to $U_m$ was evaluated. The derivative results at some $U_m$ values are given in Fig. 16. Each boxplot in this figure represents the derivative values for all 5000 points at a fixed α and β (values given on the horizontal axis). We can notice little variation for each boxplot, which means that the second derivatives are very similar all over the $X_g$ search space. Moreover, the derivative values are all close to zero. This confirms that the mean torque can be considered as a linear function with respect to $U_m$:

$$f_1(X_g, U_m) \approx g(X_g) + A(X_g)U_m \quad (8)$$

Having shown that the metamodel constructed for the mean torque can be approximated to a linear function with respect to $U_m$, let us calculate the variability of the operator $A(X_g)$ with respect to $X_g$. According to our linear approximation, the operator $A(X_g)$ is equivalent to the partial derivative with respect to $U_m$ of the mean torque's metamodel. Fig. 17 shows the results in the same manner as Fig. 16, this time for the first derivative with respect to $U_m$. Once again, no significative variations are observed for all the boxplots meaning that the first partial derivative with respect to $U_m$ is almost independent of $X_g$. The matrix $A(X_g)$ can be considered as independent of $X_g$:

$$f_1(X_g, U_m) \approx g(X_g) + AU_m \quad (9)$$

In the same way, we observe that this first derivative is also almost independent of $U_m$ (the derivative values are very close for all the boxplots). This is consistent with what was mentioned about the second derivative in Fig. 16.

So far, we have shown the linear behavior of the mean torque ($f_1$) with respect to $U_m$. However, this doesn't yet explain the results obtained in Fig. 13 and Fig. 14. According to our approximation, the expectation of $f_1$ with respect to $U_m$ and $U_g$ as well as its expectation with respect to $U_g$ alone are written as follows:

$$E_{(Ug,Um)}[f_1(X_g+U_g,U_m)] \approx E_{Ug}[g(X_g+U_g)] + A\mu_{Um} \quad (10)$$

$$E_{Ug}[f_1(X_g+U_g,U_{nom})] \approx E_{Ug}[g(X_g+U_g)] + AU_{nom} \quad (11)$$

with $\mu_{Um}$ the mean of $U_m$ uncertain variables. And, for the worst-case:

$$WC_{(Ug,Um)}[f_1(X_g+U_g,U_m)] \approx WC_{Ug}[g(X_g+U_g)] + WC_{Um}[AU_m] \quad (12)$$

$$WC_{Ug}[f_1(X_g+U_g,U_{nom})] \approx WC_{Ug}[g(X_g+U_g)] + AU_{nom} \quad (13)$$

We note that the approximate expression for the expectation of $f_1$ with respect to $U_g$ (11) where $U_m$ is equal to the nominal material properties $U_{nom}$ and the same expectation, but with respect to $U_g$ and $U_m$ (10) are the same except for terms that do not depend on $X_g$. The same is true for the worst-case expressions of $f_1$ (12) and (13). This means that their maximum values will be for the same $X_g$ but with shifted torque values (for the case of $U_g$ and $U_m$ with respect to the case of $U_g$ only).

*3.3  Results verification by FEM simulation*

All the results presented in the previous sections were based on surrogate models. We will therefore verify in this section some results using FEM simulations. Fig. 18 and Fig. 19 show the same boxplots based on the metamodels as in Fig. 11 and Fig. 12. However, this time, FEM simulations were used to produce the same boxplot for comparison. Only two zones are represented in these figures : $A_{Exp}^{Ug}$ and $E_{Exp}^{Ug}$. The STD values were also added for compaison.

The results from Fig. 18 confirm what we have already noticed using the meta-models. For both zones $A_{Exp}^{Ug}$ and $E_{Exp}^{Ug}$, the robust optimization presents similar or better solutions than the deterministic optimization in terms of robustness. For zone $A_{Exp}^{Ug}$, both machines have practically the same FEM-computed mean torque ripple with a smaller STD for the robust design. As for zone $E_{Exp}^{Ug}$, the robust design presents a lower expected mean torque. The FEM-computed q3 value of the robust design is also equal to the median of the deterministic one (5.8%), meaning that 75% of the produced machines using the robust design would have a torque ripple lower than 5.8% vs. only 50% for the deterministic design. Although the results in Fig. 18 were obtained from the robust optimization based on the expectance formulation, the design in zone $E_{Exp}^{Ug}$ also present a particularly lower torque ripple in a worst-case scenario, with 8.4% vs. 9.9% for the deterministic design.

As seen in section 3.1, using the optimization formulation with the worst-case scenario tends to reduce the STD of the robust designs compared to the deterministic ones, especially at low torque. This is also confirmed by FEM simulations in Fig. 19. For the machines in zone $A_{WC}^{Ug}$, the robust design has a notably lower FEM-computed STD compared to the deterministic one (0.4% vs 0.78%). As for the worst-case torque ripple value, both designs possess similar FEM-computed performances, around 7.3%. As for zone $E_{WC}^{Ug}$, similar results were obtained with deterministic and worst-case optimizations with a torque ripple of 11.8% and 11.6% respectively.

Regarding the precision of the meta-models compared to FEM simulations, some differences can be noticed. This lack of precision, especially for torque ripple, is somehow expected with a strategy using fixed metamodels in optimization, especially with the difficulty to fit a torque ripple metamodel as seen in section 2.1. A vast number of simulations is needed in this case to obtain an acceptable level of accuracy. An alternative approach could be to use an adaptative strategy to update the surrogate model with additional simulations during the optimization. The additional simulations will guide the algorithm towards the optimal zones while improving the meta-models' accuracy. By using this kind of approach, additional FEM simulations are performed only in promising points in the search space, limiting computational time and increasing the precision for optimal designs. This will be the subject of a future study.

Finally, we note that surrogate-based optimizations require much less computational time compared to FEM-based optimizations. For example, 17 hours was needed to finalize the DOE simulations *(Intel(R) Xeon(R) W-2195 CPU @ 2.30 GHz and 18 cores)*. The deterministic optimization using the meta-model took only 13 seconds for 300 iterations and 150 particles. Doing the same optimization using FEM simulations would have taken around 3 months to complete.

## 4  CONCLUSIONS

We have presented in this paper a comparison between different optimizations performed on a Permanent Magnet assisted Synchronous Reluctance Machine. The first optimization used a deterministic formulation considering certain all the design parameters. In order to study the sensitivity of such design to parameters uncertainties, robust optimizations were performed considering some geometric and magnetic parameters as uncertain. Two different formulations were adopted for the robust optimizations. In the first one, we used the expectations of the mean torque and torque ripple as objective functions. In the second formulation, the objective functions were defined as the worst-case values. In order to reduce computation time, surrogate models have been built for each of the objective functions. These surrogate models have been also used to perform a Sobol' indices-based sensitivity analysis to detect the most impacting input parameters. Objective functions' expectations were computed with a quasi-Monte Carlo scheme while worst-cases where calculated with an embedded Particle Swarm Optimization algorithm. Both geometrical and magnetic property tolerances were shown to substantially influence the machines' performances, with the magnetic properties tolerances having a greater impact on the mean torque. More precisely, we found that the magnetic material properties tolerances affect the mean torque linearly. We have shown than the considered magnetic material degradation has almost the same impact on the performances over all the search space: The Pareto front of the robust optimization considering both type of uncertainties (geometrical $U_g$ and magnetic material properties $U_m$) is simply obtained by shifing horizontally (toward lower torque values) the Pareto front of the robust optimization with only the geometrical parameters ($U_g$) as uncertain. This result cannot be generalized since the search space in this paper has been restricted to a small region around the reference design.

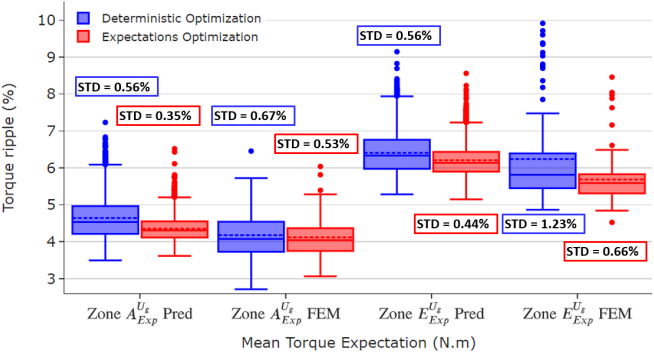

Fig. 18: Boxplots showing comparisons of torque ripple between FEM simulations and surrogate model predictions: deterministic (blue) and robust (red) machines with similar predicted Mean Torque/ Mean Torque Expectation values.

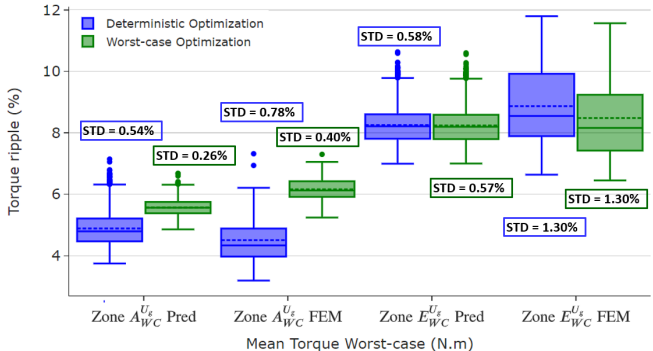

Fig. 19: Boxplots showing comparisons of torque ripple between FEM simulations and surrogate model predictions: deterministic (blue) and robust (green) machines with similar predicted Mean Torque Worst-Case values.

It should also be noted that while the predicted values of mean torque were consistent with Finite Element Method simulations, some differences were observed for torque ripple. This problem will be addressed in future work by using an adaptive strategy, like Bayesian Optimization, to update the surrogate models with additional simulations during the optimization. Nevertheless, the Finite Element Method simulations have globally confirmed the trend of the predicted results using the metamodels. The comparison of Pareto fronts has shown that robust solutions outperform deterministic solutions in terms of different robustness criteria as quantiles, expectation, worst-case and standard deviation values. However, in a constrained optimization problem, it is recommended to use the worst-case as a constraint and not as an objective. In a real problem like a production factory, the worst-case, which is very conservative, is replaced by a quantile. Such formulation allows increasing the average performance of a batch of prototypes for example, while ensuring that only a small percentage of them would not respect the defined requirements.

The aforementioned aspects show the importance of developing and applying new techniques of optimization of electrical machines when dealing with uncertainties.

## ACKNOWLEDGMENTS

We would like to thank the CONAHCYT (Consejo Nacional de Humanidades Ciencias y Tecnologías) for its support.